\documentclass[11pt,a4paper]{article}

\usepackage{amsmath}
\usepackage{amssymb}
\usepackage{amsthm}
\usepackage{amsfonts}
\usepackage{mathrsfs}

\newcommand{\N}{\mathbb{N}}
\newcommand{\mfrak}{{\mathfrak m}}

\newcommand{\codim}[0]{\operatorname{codim}}

\newcommand{\reg}[0]{\operatorname{reg}}

\newcommand{\supp}[0]{\operatorname{supp}}
\newcommand{\udeg}[0]{\operatorname{\underline{deg}}}

\newtheorem{satz}{Theorem}[section]
\newtheorem{lemma}[satz]{Lemma}

\newtheorem{notation}[satz]{Notation}
\newtheorem{coro}[satz]{Corollary}
\newtheorem{example}[satz]{Example}

\newtheorem{remarks}[satz]{Remarks}

\newtheorem*{frago}{Question}

\newtheorem*{definitiono}{Definition}

\title{Castelnuovo-Mumford regularity and reduction number of smooth monomial curves}

\author{Michael Hellus$^1$,  L\^{e} Tu\^{a}n Hoa$^{2,}$\footnote{The second named  
author was supported by the National Basic Research Program
(Vietnam) and Max-Planck Institute for Mathematics in the Sciences (Germany). He
would like to thank the MIS for the financial support and the hospitality.}, J\"urgen 
St\"uckrad$^3$\\ {\small $^{1,3}$Universit\"at Leipzig, Fakult\"at f\"ur Mathematik und 
Informatik},\\
{\small Augustusplatz 10/11, D-04109 Leipzig, Germany }\\ {\small $^2$Institute of 
Mathematics Hanoi, }\\ {\small 18 Hoang Quoc Viet Road, 10307 Hanoi, Vietnam}\\
{\small E-mail:  $^1$Michael.Hellus@math.uni-leipzig.de,  $^2$lthoa@math.ac.vn,}\\ 
{\small $^3$stueckrad@math.uni-leipzig.de}}

\begin{document}
\maketitle

\begin{abstract}
We compare, for smooth monomial projective curves, the
Castel\-nuo\-vo-Mumford regularity and the reduction number; we
pre\-sent an example where these two numbers differ. However, we
show they coincide for a certain class of monomial curves.
Furthermore, for smooth monomial curves we prove an inequality which
is stronger than the one from the Eisenbud-Goto conjecture.
\end{abstract}
\section{Introduction}
The Eisenbud-Goto conjecture states that
\[ \reg R\leq \deg R-\codim R\]
holds for every graded domain $R$ over an algebraically closed base
field $K$, where $\reg R$ is the smallest integer $n$ such that
$H^i_\mfrak(R)_{j-i}=0$ for all $i\geq 0$ and $j>n$, $\deg R$ is the
multiplicity of $R$ and $\codim R:=\dim _K (R_1)-\dim R$. This
conjecture was proved in various special cases; \cite{gruson83}
contains a proof for projective irreducible curves.

On the other hand, \cite{hoa03} contains a proof of $r(S)\leq \deg
K[S] -\codim K[S]$ where $S$ is a (homogenous) simplicial affine
semigroup and $r(S)$ is the reduction number of $K[S]$ with respect
to the natural minimal reduction of the simplicial ring $K[S]$.
Moreover, one has a rather close relation between $r(S)$ and $\reg K[S]$, namely $r(S)\leq \reg 
K[S] < (\dim K[S]) r(S)$ (by \cite{trung87} and \cite{hoa03}). Therefore, it is natural to ask

\begin{frago}

Does \[ r(S)=\reg K[S]\tag{$Q$}\] hold?

\end{frago}

Besides of own interest, a positive answer to this question would confirm the Eisenbud-Goto 
conjecture for $\reg K[S]$.  As this is maybe the simplest case, we investigate $(Q)$ for smooth 
monomial curves. We get partial positive answers (Theorem
\ref{partial}, see also Remark \ref{sechzehn} (ii)). Unfortunately, in
general the answer is negative (see Example \ref{counterexample}). In order to get these 
results, we have to reformulate the above problem to some combinatorial problems. Then we 
can, in particular, use the software {\it Macaulay 2} to search for a counterexample.

Though $(Q)$ has negative answer in general, we can establish a new  upper
bound (Theorem \ref{improvement}) for the Castelnuovo-Mumford regularity $\reg K[S]$ of 
smooth
monomial curves that is much  stronger than what the Eisenbud-Goto
conjecture claims. Note that \cite{herzog03} contains a
combinatorial proof of the Eisenbud-Goto conjecture  for simplicial
toric rings with isolated singularity.

Finally, in some cases we can compute the reduction number $r(S)$
(see Theorem \ref{computeR}).

\section{Combinatorial formulation}

We will keep the following notation for the rest of this article. Let $\N,\ \mathbb{Z}_+$ and
$\mathbb{Z}_-$ denote the set of non-negative, positive and negative integers,
respectively.
Let $S\subset \N^2$ be a  semigroup generated by
\[ {\cal A}=\{ {\bf e}_1 := (\alpha, 0),\ {\bf e}_2= (0,\alpha), \ (\alpha - a_1, a_1),\dots , 
(\alpha-a_c,a_c)\} ,\]
where $c\geq 1,\alpha \geq 2$ and $0 < a_1 <\cdots < a_c< \alpha$ are relatively prime 
integers. Note that $c=\codim
K[S]$. The degree of an element of the group generated by $S$ is
defined as the sum of its entries, divided by $\alpha $. By $K[S] \equiv K[t_1^\alpha, 
t_1^{\alpha-a_1}t_2^{a_1},..., t_1^{\alpha-a_c}t_2^{a_c}, t_2^\alpha] \subset  K[t_1,t_2]$  
we denote the affine semigroup ring
associated to $S$. Then ${\mathfrak q} := (t_1^\alpha,t_2^\alpha)$ is a minimal reduction of 
the maximal homogeneous ideal $\mfrak$ of $K[S]$. Let $r(S)$ denotes the reduction number 
of $\mfrak$ with respect to $\mathfrak q$, that is $r(S) = \min\{ r\in \N|\ \mfrak^{r+1} = 
{\mathfrak q}\mfrak^r\}$.

\begin{remarks}

\label{bem}{\rm (i) Note that one can easily compute $r(S)$ as follows:
 \[ r(S)=\min \{ r\in \N  \vert (r+1){\cal A}=\{ {\bf
e}_1, ,{\bf e}_2\} +r{\cal A}\} \geq 1.\]
(Notation:  ${\cal B}+{\cal C}:=\{ b+c\vert \ b\in {\cal B},c\in
{\cal C}\}$ for subsets $\cal B,C$ of $\N^2$).

(ii) One has
$ r(S)\leq \reg K[S] $
(see \cite{trung87}) and
\[ \reg K[S]\begin{cases}\leq 2r(S)-2,  & \mbox{ if } r(S)\geq 2\\=1,  & \mbox{ if 
}r(S)=1\end{cases}\]
(see \cite[Theorem 3.1]{hoa03}). Thus, if $r(S)\leq 2$ or $\reg
K[S]\leq 3$, then $r(S)=\reg K[S]$. The ``smallest" counterexample to
$(Q)$ could  therefore only occur  for $r(S)=3$ and $\reg K(S)=4$.

(iii) For every element $b=(b_1,b_2)$ of $\cal A$ one has
$b_1+b_2=\alpha $; hence we can simplify our notation by replacing
every element of $\cal A$ by its first entry, that means
\[ {\cal A} = \{ 0, \alpha- a_c,...,\alpha-a_1, \alpha\}.\]
\hfill $\square $}

\end{remarks}

\begin{definitiono}

(i) We say a subset ${\cal B}$ of $\N$ is {\rm full} iff it has the
form
\[ {\cal B}=\{ 0,1,\dots ,m\} \] for some $m\in \N$.

(ii) We say a subset $\cal A$ of $\N$ {\rm has} $(P_1)$ iff for
every $m\in \N$ one has
\[ m{\cal A}\text{ not full}\Rightarrow m{\cal
A}+\{ 0,\alpha \} \text{ not full}. \]

(iii) We say a subset $\cal A$ of $\N$ {\rm has} $(P_2)$ iff for
every $m\in \N$ one has
\[ m{\cal A}\text{ not full}\Rightarrow   m{\cal
A}\text{ does not contain all numbers from the list }\]
\[ 0,1,\dots ,\alpha ,(m-1)\alpha , (m-1)\alpha +1,\dots ,m\alpha .\]

\end{definitiono}

If ${\cal A}=\{ 0,p_1,\dots ,p_l,\alpha \} $, $0<p_1<\dots
<p_l<\alpha $ and $m\alpha $ is full for some $m\geq 1$ then
$1,\alpha -1\in {\cal A}$; this means that $\cal A$ defines a smooth
monomial curve. Also, if $m{\cal A}$ is full then $n{\cal A}$ is
full for every $n>m\geq 1$.

For the rest of the section we assume that $\cal A$ defines a smooth
monomial curve, i.e. ${\cal A} =\{ 0, 1, p_1,...,p_l, \alpha -1,
\alpha \}$ and $1<p_1<\dots <p_l<\alpha -1$. For later use we set ${\bf g}_i:=(\alpha -i,i)$, 
where  $0\leq i \leq \alpha$. Thus, ${\bf e}_1 = {\bf g}_0$ and ${\bf e}_2 = {\bf g}_\alpha$.

\begin{lemma}\label{lemmaA} Let $G\subseteq \mathbb{Z}^2$ be the group
generated by all $(u_1,u_2)$ such that $u_1+ u_2 = \alpha $. Let
$\bar{S} = G\cap \N ^2$. Then as $\mathbb{Z}$-graded modules we have

(i) $H^2_\mfrak(K[S]) \cong K[G\cap \mathbb{Z}_-^2]$.

(ii) $H^1_\mfrak (K[S]) \cong K[\bar{S}\setminus  S]$.
\end{lemma}
{\it Proof.} Clearly ${\bf u}  := (u_1,u_2)\in S-\N {\bf e}_1$
implies $u_2 \geq 0$. Conversely, let ${\bf u} \in G$ such that $u_2
\geq 0$. Note that $G$ is generated by ${\bf e}_1, {\bf g}_1$.
Writing ${\bf u} = p{\bf e}_1 + q{\bf g}_1$, $p,q\in \mathbb{Z}$,
and comparing the second coordinates we get that $q\geq 0$. Hence
${\bf u} \in S- \N {\bf e}_1$. Thus \[ S-\N {\bf e}_1 = \{ {\bf u}
\in G| \ u_2 \geq 0\}.\] Similarly, \[ S-\N {\bf e}_2 = \{ {\bf u}
\in G| \ u_1 \geq 0\}.\] By \cite[Corollary 3.8]{trung86} we get \[
H^2_\mfrak(K[S]) \cong K[G \setminus ((S-\N {\bf e}_1)\cup (S-\N {\bf
e}_2))] =K[G\cap \mathbb{Z}_-^2].\] Since $\alpha S \subset \N {\bf
e}_1 + \N {\bf e}_2$, the set of ${\bf u}\in G$ such that ${\bf u} +
p(S\setminus \{ 0\}) \subseteq S $ for $p\gg 0$ is exactly the set
$(S-\N {\bf e}_1)\cap (S-\N {\bf e}_2) = \bar{S}$. By
\cite[Corollary 3.4(ii)]{trung86} get (ii). \hfill $\square$
\vskip0.3cm

Note that $\bar{S}$ coincides with the so-called normalization of
$S$ (see \cite{hochster72}). From Lemma \ref{lemmaA}(ii) we see that
$K[S]$ is a Cohen-Macaulay ring if and only if $\cal A$ is full. In
this case
\[ r(S) = \reg K[S] = 1.\] The first equality follows from
\cite[Proposition 3.2]{trung87}, while the second one follows from
Lemma \ref{lemmaA}.

\begin{coro}\label{coroB}We have
\[ \begin{array}{ll} \reg K[S] &= \max\{\deg {\bf u}|\ {\bf u} \in (\bar{S}\setminus S)\cup 
\{0\}\} +1 \\
 &= \min\{ m>0|\ m{\cal A}\ is\ full\}.
\end{array}\]
 \end{coro}

\begin{satz}

\label{interpr}
\[ {\cal A}\text { has }(P_2)\Rightarrow {\cal A}\text { has
}(P_1)\iff (Q)\text{ has  positive answer.}
\]

\end{satz}

{\it Proof.} First implication: If $m\cal A$ does not contain a
number $i$ from the list $0,1,\dots ,\alpha $, then $m{\cal A}+\{
0,\alpha \}$ does not contain the number $i$, too; if $m\cal A$ does
not contain a  number $i$ from the list $(m-1)\alpha , (m-1)\alpha
+1,\dots ,m\alpha $, then $m{\cal A}+\{ 0,\alpha \}$ does not
contain the number $i+\alpha $.

${\cal A}$ has ($P_1$)  $\Longrightarrow $ (Q) has positive answer:
Let $m = \reg K[S] -1$. By Corollary \ref{coroB}, $m{\cal A}$ is not
full, but $(m+1){\cal A}$ is full. If $m\geq r(S)$, then by the
definition of $r(S)$ we have $(m+1){\cal A}= m{\cal A}+\{ 0, \alpha
\}$ and $m{\cal A}+\{ 0, \alpha \}$ is full. This contradicts the
property ($P_1$) of ${\cal A}$. Hence $m\leq r(S)-1$ or equivalently
$\reg K[S] \leq r(S)$, which forces $\reg K[S] = r(S)$.

(Q) has positive answer $\Longrightarrow $ ${\cal A}$ has ($P_1$):
Assume that $m{\cal A}$ is not full. If $m{\cal A}+\{0,\alpha \}$ is
full, then $(m+1){\cal A}= m{\cal A}+\{0,\alpha \}$. This means
$r(S) \leq m$. Since $\reg(S) = r(S) \leq m$, by Corollary
\ref{coroB}, $m{\cal A}$ is full, a contradiction.  \hfill $\square$

\section{Results}
We keep all notation from the previous section.
\begin{notation}
{\rm  Set $H := (\{0,1\},+,\cdot)$, $0+0 :=0$, $0+1 = 1+0 = 1+1 :=
1$ and $0\cdot 0 := 0 \cdot 1 := 1 \cdot 0 := 0$, $1\cdot 1 := 1$.
$H$ is a commutative semiring with identity. Let $t$ be an
indeterminate. Then $H[t]$ is also a commutative semiring with
identity. Note that, for $f, g \in H[t]$ and $m \in{\mathbb Z}_+$, one has
$(f+g)^m = f^m + f^{m-1}g + \cdots + g^m$. Furthermore one has
$\supp (f+g) = \supp f \cup \supp g$. For $f \ne 0$ set $\deg f :=
\max\supp f$ and $\udeg f := \min\supp f$.

For $p \in \N$ we define $h_p := \sum_{i=0}^p t^i$. Then $(h_p)^i =
h_{ip}$ holds for all $i \in \N$.

Let $f \in H[t]$. By definition, a {\it  gap} of $f$ is a subset $L$
of $\N \setminus\supp f$ that has the form $L = \{i,i+1,\dots,j\}$
where $i, j \in \N$, $i \le j$ and such that $i-1,j+1 \in \supp f$.
If $L$ is a gap of $f$, its {\rm length} is $\sharp L$, by
definition. We set
\[\lambda(f) := \max\{\sharp L~|~L \mbox{ gap of } f\}.\]
}
\end{notation}

\begin{lemma} \label{lemmajuergen}For $f, g \in H[t]$ the following statements hold:

(i) $\lambda(h_pf) = \max\{0,\lambda(f) - p\}$ for all $p \in \N$.

(ii) $\lambda(f+g) \le \max\{\lambda(f),\lambda(g),\udeg f - \deg g
- 1,\udeg  g - \deg f - 1\}$.

(iii) $\lambda(fg) \le \max\{\lambda(f),\lambda(g)\}$.

(iv) For every $i \in {\mathbb Z}_+$ one has $\lambda(f^i) \le
\lambda(f)$.

\end{lemma}

For a subset $\cal A$ of $\N$ we define
\[f_{\cal A} := \sum_{i = 0}^{\alpha} \epsilon_i t^i \in H[t],~\mbox{ where }~\epsilon_i := 
\begin{cases}1,  & \mbox{if } i \in {\cal A}.\\0,  & \mbox{else.}\end{cases}\]
Then $f_{m\cdot \cal A} = f_{\cal A}^m$ holds for every $m\in \N$.

\begin{lemma}

\label{komb} $(P_1)$ does not hold in general; e.~g., $(P_1)$ does
not hold for $m=3$ and
$ {\cal A}:=\{ 0,1,2,5,13,14,16,17\}. $

(ii) $(P_2)$ holds if $\cal A$ has the form
\[ \{ 0,1,\dots ,p,q=q_1,\dots ,q_l,\alpha \} \]
with $1\leq p<q_1<\dots <q_l\leq \alpha $ and $p\geq \alpha -q$. (By
symmetry, $(P_2)$ then also holds if $\cal A$ has the form
\[ \{ 0,q_1,\dots ,q_l=q,\alpha -p,\alpha -p+1,\dots ,\alpha \} \]
with $0<q_1<\dots <q_l<\alpha -p\leq \alpha -1$ and $p\geq q$).

\end{lemma}

{\it Proof.} (i) One calculates
$ 3{\cal A} =\{ 0,1,\dots ,51\} \setminus \{ 25\}$
and $ 3{\cal A}+\{ 0,17\} =\{ 0,1,\dots ,68\}$.  This example
was found with the software {\it Macaulay 2} (\cite{macaulay2}).
\smallskip

(ii) $f := f_{\cal A}$ has the form $f = h_p + t^q g$ with $g \in
H[t]$ and $\deg g = \alpha - q \le p$ and $g(0) = 1$.

For $m \in {\mathbb Z}_+$ one has $f^m = h_p^m + t^qh_p^{m-1}g + \cdots +
t^{mq}g^m = h_{mp} + t^qh_{(m-1)p}g + \cdots + t^{mq}g^m$. Let $i
\in \{0,\dots,m-1\}$. Then from Lemma \ref{lemmajuergen} we deduce
that
\begin{eqnarray*}
\lambda(t^{iq}h_{(m-i)p}g^i) &=& \lambda(h_{(m-i)p}g^i) = \max\{0,\lambda(g^i) - 
(m-i)p\}\\
&\le& \max\{0,\lambda(g) - p\} \le \max\{0,\deg g - p\}\\
&=& 0,
\end{eqnarray*}
i.~e. $t^{iq}h_{(m-i)p}g^i = t^{iq}h_{\beta_i}$ with $\beta_i := mp
- i(p+q-\alpha)$. We have
\begin{eqnarray*}
\udeg(t^{(i+1)q}h_{(m-i-1)p}g^{i+1}) - \deg (t^{iq}h_{(m-i)p}g^i) &=& (i+1)q - iq - \beta_i\\
&=& i(p+q -\alpha) + q -mp.
\end{eqnarray*}
Assume $\lambda(f^m) > 0$. From what was shown above it follows that
$i(p+q -\alpha) + q -mp
> 1$ for some $i \in \{0,\dots,m-1\}$ or $\lambda(g^m) > 0$, i.~e. precisely one of the 
following two conditions holds:
\begin{enumerate}
\item[(A)] $(m-1)(p+q -\alpha) + q -mp > 1,$
\item[(B)] $(m-1)(p+q -\alpha) + q -mp \le 1$ and $\lambda(g^m) >0$.
\end{enumerate}
Case (A): Then $mq-1 \not \in \supp f^m$. Because of $mq -
(m-1)\alpha -p = (m-1)(p+q -\alpha) + q -mp > 1$ we have $mq - 1 >
(m-1)\alpha + p \ge (m-1)\alpha$ and, therefore, ($P_2$) holds.

Case (B): There exists $j\in \N$ such that
$(m-1)\alpha + p = \deg
(t^{(m-1)q}h_p g^{m-1}) $ $ < j \le \deg f^m$ and $j \not\in \supp f^m$.
Because of $j > (m-1)\alpha + p \ge (m-1)\alpha$ property ($P_2$)
holds in this case, too.\hfill $\square $

\begin{example}
\label{counterexample} In general $(Q)$ has negative answer; e.~g.,
over an arbitrary field $K$, let $S$ be the subsemigroup of $\N^2$
generated by
\[ \{ (17,0),(16,1),(1,16),(0,17),(14,3),(13,4),(5,12),(2,15)\} .\]
Then
\[ r(S)=3<4=\reg K[S] \]
This follows from Lemma \ref{komb}(i) and Theorem
\ref{interpr}.\hfill $\square $
\end{example}

\begin{remarks} {\rm
\label{sechzehn} (i) The preceding example is the "smallest"
possible counterexample to $(Q)$ in the sense of Remark
\ref{bem}(ii).

(ii) Using {\rm Macaulay 2} one can check that $(P_2)$ holds for
$m=3$ and for every set
\[ {\cal A} =\{ 0, 1, p_1,...,p_l, \alpha -1,
\alpha \} \] ($1<p_1<\dots <p_l<\alpha -1$) and $\alpha \leq 16$.
This has the following geometric meaning: If, over an arbitrary
field $K$, the subsemigroup $S$ of $\N^2$ corresponds to a smooth
(monomial) curve of degree at most 16, then the following property
holds
\[ \reg K[S]\geq 4\Rightarrow r(S)\geq 4. \] }
\hfill $\square $
\end{remarks}

\begin{satz}
\label{partial} Over an arbitrary field $K$, let $S\subseteq \N ^2$
be the subsemigroup associated to
\[ {\cal A}=\{ 0,1,\dots ,p,q=q_1,\dots ,q_l,\alpha \}
\] with $1\leq p<q_1<\dots <q_l= \alpha -1$ and $p\geq \alpha -q$.
Then
$ r(S)=\reg K[S]$.

\end{satz}
{\it Proof. }This follows from Theorem \ref{interpr} and Lemma
\ref{komb}(ii).\hfill $\square $ \vskip .2cm

By a famous result of Gruson, Lazarsfeld and Peskine
\cite{gruson83}, the Castelnuovo-Mumford regularity of a projective
irreducible curve $C$ is bounded by $\deg C - \codim C +1$. For a
smooth monomial curve defined by ${\cal A}$ this means that $\reg
K[S] \leq (\sum \sharp L) + 1,$ where $L$ runs over all gaps of
$f_{{\cal A}}$. In \cite[Corollary 2.2]{herzog03} there is a
combinatorial of this result in this case. We can improve this bound
as follows:

\begin{satz}
\label{improvement} Let $\varepsilon = \max\{ i|\ 0,1,...,i,\ \alpha
, \alpha -1,...,\alpha -i\in {\cal A}\} \geq 1$. Then \[ \reg K[S]
\leq \left[\frac{\lambda (f_{{\cal A}}) -1}{\varepsilon }\right]
+2,\] where $[a]$ denotes the integer part of $a\in \mathbb{R}$.
\end{satz}

 {\it
Proof. }If $K[S]$ is a Cohen-Macaulay ring, then by Lemma
\ref{lemmaA}, $\lambda (f_{\cal A}) =0$ and the claim follows by
Corollary \ref{coroB}.

Assume  $\bar{S} \neq S$. Recall that  ${\bf g}_i = (\alpha -i, i)$.  Let $u \in
\bar{S} \setminus S$. Then we can write ${\bf u} = {\bf y} + p{\bf
e}_1 + q{\bf e}_2$, where $p,q\in \N ,\ {\bf y}  = (y_1,y_2)\in \N
^2$ and $y_1 + y_2 = \alpha $. Since ${\bf u} \not\in S,\ {\bf y} +
i{\bf e}_1 \not\in S$ for all $i= 0,...,p$. Note that \[ {\bf y} +
i{\bf e}_1 = (i-j) {\bf g}_t + j {\bf g}_{t+1} + (y_1+ti+j, y_2
-ti-j),\] where $0\leq j \leq i \leq p$ and $0\leq t\leq \varepsilon
-1$. By the definition of $\varepsilon $ we have ${\bf g}_0,...,{\bf
g}_\varepsilon  \in S$. Hence  for all $i,j,t$ as above we must have
$(y_1+ti+j, y_2 -ti-j)\not\in S$. This means all numbers $y_1,
y_1+1,...,y_1+\varepsilon p$ belong to a gap $L$ of $f_{{\cal A}}$.
Similarly, the condition $y+ q{\bf e}_2 \not\in S$ forces
$y_1,y_1-1, ...., y_1-\varepsilon q \in L$. Therefore $\sharp L \geq
1 + (p+q) \varepsilon $. Since $\deg {\bf u} = p+q+1$,  we get
\[ \deg {\bf u}   \leq \left[\frac{\sharp L
-1}{\varepsilon }\right] + 1  \leq \left[\frac{\lambda (f_{{\cal A}})
-1}{\varepsilon }\right] +1 .\]
By Corollary \ref{coroB}, the claim follows. \hfill $\square$

The following result gives a lower bound  for the reduction number
and the Castelnuovo-Mumford regularity in terms of the length of the
first gap of $\lambda (f_{{\cal A}})$.  It gives further cases when (Q) holds.  We also see  how 
one can in
some cases compute the reduction number $r(S)$ using the information
on $\reg K[S]$.

\begin{satz}
\label{computeR}
 Assume that \[ {\cal A} =\{ 0,1,..., \varepsilon ,
p_1,...,p_l , \alpha - 1,\alpha \},\] such that $\varepsilon +2 \le p_1 <\cdots < p_l
\leq \alpha - \varepsilon$. Then \[ \reg K[S] \geq r(S) \geq
\left[\frac{p_1 -2}{\varepsilon }\right] + 1.\]
In particular, if also $\alpha - i\in {\cal A}$ for all $0\leq i\leq \varepsilon$ and
$p_1- \varepsilon -1 = \lambda (f_{{\cal A}})$, then
\[ \reg K[S] = r(S) = \left[\frac{\lambda (f_{{\cal A}})
-1}{\varepsilon }\right]+ 2.\]
\end{satz}
{\it Proof. }
 Let
\[  \delta = \left[\frac{p_1-\varepsilon -2}{\varepsilon }\right] \
\text{and} \ {\bf u} := (u_1,u_2) = (p_1-1, \alpha -p_1 +1) + \delta
(0, \alpha ).\] If ${\bf u} \in S$, then ${\bf u} = \sum_{0\leq
i\leq \alpha} q_i(i,\alpha -i)$, where $\sum q_i = \delta +1$.
Comparing the first coordinates we get $q_i= 0$ for all
$i>\varepsilon $. Hence \[ u_1  = \sum_{i=0}^{\varepsilon } i q_i
\leq \varepsilon \sum_{i=0}^{\varepsilon }  q_i = \varepsilon
(\delta +1) \leq   p_1-2 < p_1-1, \]
a contradiction. Hence ${\bf u} \not\in S$. On the other hand,
letting $p_1-\varepsilon -1= \delta \varepsilon + \gamma $, where
$0\leq \gamma <\varepsilon $ we get
 \[ \begin{array}{ll}
{\bf u} + {\bf e}_2  & = (p_1-1, \alpha -p_1-1) + (\delta +1)(0,\alpha ) \\
&= ((\delta +1)\varepsilon + \gamma , \alpha + (\delta +1)(\alpha -\varepsilon ) - \gamma )\\
&= (\delta +1)(\epsilon , \alpha - \varepsilon ) + (\gamma , \alpha
- \gamma ) \in S.
\end{array}\]
 Since ${\bf u} \not\in S$, ${\bf u} + {\bf e}_2 \not\in {\bf e}_2+S$. Comparing the first 
coordinate we also get ${\bf u} + {\bf e}_2 \not\in {\bf e}_1 + S$. Thus ${\bf u} + {\bf e}_2 
\not\in \{{\bf e}_1, {\bf e}_2\} + S$ and $r(S) \geq \deg ({\bf u} + {\bf e}_2) = \delta +2$. 
Since we always have $r(S) \leq \reg K[S]$, this proves the first statement.

If $p_1- \varepsilon -1 = \lambda (f_{{\cal A}})$, then  combining
with Theorem \ref{improvement}, we finally get \[  \delta +2 \leq r(S) \leq  \reg
K[S] \leq   \left[\frac{\lambda (f_{{\cal A}}) -1}{\varepsilon
}\right] +2=   \left[\frac{p_1-\varepsilon -2}{\varepsilon }\right] + 2 = \delta +2,\] which yields 
\[  r(S) = \reg K[S] = \delta
+2 =  \left[\frac{\lambda (f_{{\cal A}}) -1}{\varepsilon }\right]
+2.\] \hfill $\square$

\end{document}